\documentclass[11pt]{amsart}
\usepackage{graphicx}

\newtheorem{theorem}{Theorem}
\newtheorem{lemma}[theorem]{Lemma}
\newtheorem{proposition}[theorem]{Proposition}

\numberwithin{equation}{section}

\begin{document}

\title{On a Conjecture of Milnor about Volume of Simplexes}

\author{Ren Guo}

\address{Department of Mathematics, Rutgers University, Piscataway, NJ, 08854}

\email{renguo@math.rutgers.edu}

\author{Feng Luo}

\address{The Center of Mathematical Science, Zhejiang
University, Hangzhou, China, 310027}

\address{Department of Mathematics, Rutgers University, Piscataway, NJ, 08854}

\email{fluo@math.rutgers.edu}


\date{}

\dedicatory{Dedicated to the memory of Xiao-Song Lin.}


\begin{abstract} We establish the second part of Milnor's
conjecture on the volume of simplexes in hyperbolic and spherical
spaces. A characterization of the closure of the space of the
angle Gram matrices of simplexes is also obtained.
\end{abstract}
\maketitle

\section{Introduction}

\subsection*{Milnor's conjecture}

In \cite{Mi}, John Milnor conjectured that the volume of a
hyperbolic or spherical $n-$simplex, considered as a function of
the dihedral angles, can be extended continuously to the
degenerated simplexes. Furthermore, he conjectured that the
extended volume function is non-zero except in the closure of the
space of Euclidean simplexes. The first part of the conjecture on
the continuous extension was established in \cite{Lu2} (\cite{Ri}
has a new proof of it which generalizes to many polytopes). The
purpose of the paper is to establish the second part of Milnor's
conjecture.

To state the result, let us begin with some notations and
definitions. Given an $n-$simplex in a spherical, hyperbolic or
Euclidean space with vertices $u_1, ..., u_{n+1}$, the $i$-th
codimension-1 face is defined to be the $(n-1)$-simplex with
vertices $u_1,...,u_{i-1}$,    $u_{i+1},..., u_{n+1}$. The
dihedral angle between the $i$-th and $j$-th codimension-1 faces
is denoted by $\theta_{ij}$. As a convention, we define
$\theta_{ii} =\pi$ and call the symmetric matrix $A=[-\cos
(\theta_{ij})]_{(n+1) \times (n+1)}$ the \it angle Gram matrix \rm
of the simplex.  It is well known that the angle Gram matrix
determines a hyperbolic or spherical $n$-simplex up to isometry
and Euclidean $n$-simplex up to similarity. Let
$\mathcal{X}_{n+1}, \mathcal{Y}_{n+1}, \mathcal{Z}_{n+1}$ in $
\mathbf{R}^{(n+1) \times (n+1)}$ be the subsets of $(n+1) \times
(n+1)$ symmetric matrices corresponding to the angle Gram matrices
of spherical, hyperbolic, or Euclidean $n-$simplexes respectively.

The volume of an $n$-simplex can be expressed in terms of the
angle Gram matrix by the work of Aomoto \cite{Ao}, Kneser
\cite{Kn} and Vinberg \cite{Vi}. Namely, for a spherical or
hyperbolic $n$-simplex $\sigma^n$ with angle Gram matrix $A$, the
volume $V$ is

\begin{equation} \label{1} V(A) = \mu_n^{-1} \sqrt{|\det(ad(A))|}
\int_{ \mathbf{R} ^{n+1}_{\geq 0}} e^{ -x^t ad(A)x} dx .
\end{equation}
where $ \mathbf{R}^{n+1}_{\geq 0}=\{(x_1,...,x_{n+1})|\
x_i\geq0\}$, the constant $\mu_k = \int_0^{\infty} x^k e^{-x^2}
dx$ and $ad(A)$ is the adjoint matrix of $A$. In \cite{Lu2}, it is
proved that the volume function $V: \mathcal{X}_{n+1} \cup
\mathcal{Y}_{n+1} \to \mathbf{R}$ can be extended continuously to
the closure $\overline{\mathcal{X}}_{n+1} \cup
\overline{\mathcal{Y}}_{n+1}$ in $ \mathbf{R}^{(n+1) \times
(n+1)}$. The main result of this paper, which verifies the second
part of Milnor's conjecture, is the following theorem.

\begin{theorem} The extended volume function $V$ on
the closure $\overline{\mathcal{X}}_{n+1} \cup
\overline{\mathcal{Y}}_{n+1}$ in $\mathbf{R}^{(n+1) \times (n+1)}$
vanishes at a point $A$ if and only if A is in the closure
$\overline{\mathcal{Z}}_{n+1}$.
\end{theorem}

\subsection*{A characterization of angle Gram matrices}
We will use the following conventions. Given a real matrix
$A=[a_{ij}],$ we use $A\geq 0$ to denote $a_{ij}\geq 0$ for all
$i,j$ and $A> 0$ to denote $a_{ij}> 0$ for all $i,j.$ $A^t$ is the
transpose of $A.$ We use $ad(A)$ to denote the adjoint matrix of
$A.$ The diagonal $k\times k$ matrix with diagonal entries
$(x_1,...x_k)$ is denoted by $diag(x_1,...x_k).$ A
characterization of the angle Gram matrices in $\mathcal{X}_{n+1},
\mathcal{Y}_{n+1}$ or $\mathcal{Z}_{n+1}$ is known by the work of
\cite{Lu1} and \cite{Mi}.

\begin{proposition}[\cite{Lu1},\cite{Mi}] Given an
$(n+1)\times (n+1)$ symmetric matrix $A=[a_{ij}]$ with $a_{ii}=1$
for all $i,$ then
\begin{enumerate}
\item [(a)] $A \in \mathcal{Z}_{n+1}$ if and only if $det(A)=0$,
$ad(A)>0$ and all principal $n \times n$ submatrices of $A$ are
positive definite,

\item [(b)] $A \in \mathcal{X}_{n+1}$ if and only if $A$ is
positive definite,

\item [(c)] $A \in \mathcal{Y}_{n+1}$ if and only if $det(A) <0$,
$ad(A)>0$ and all principal $n \times n$ submatrices of $A$ are
positive definite.
\end{enumerate}
In particular, all off-diagonal entries $a_{ij}$ have absolute
values less than 1, i.e.,$|a_{ij}|<1$ for $i \neq j$.
\end{proposition}

The following gives a characterization of matrices in
$\overline{\mathcal{X}}_{n+1},\overline{\mathcal{Y}}_{n+1}$ and
$\overline{\mathcal{Z}}_{n+1}$ in $ \mathbf{R}^{(n+1) \times
(n+1)}$.

\begin{theorem} Given an $(n+1)\times (n+1)$
symmetric matrix $A=[a_{ij}]$ with $a_{ii}=1$ for all $i,$ then
\begin{enumerate}
\item [(a)] $A \in \overline{\mathcal{Z}}_{n+1}$ if and only if
$det(A)=0$, $A$ is positive semi-definite, and there exists a
principal $(k+1)\times (k+1)$ submatrix $B$ of $A$ so that
$det(B)=0, ad(B)\geq 0$ and $ad(B)\neq 0$,

\item [(b)] $A \in \overline{\mathcal{X}}_{n+1}$ if and only if
either $A$ is in $\mathcal{X}_{n+1}$ or there exists a diagonal
matrix $D=diag(\varepsilon_1,...,\varepsilon_{n+1})$ where
$\varepsilon_i=1$ or $-1$ for each $i=1,...,n+1,$ such that
$DAD\in \overline{\mathcal{Z}}_{n+1},$

\item [(c)] $A \in \overline{\mathcal{Y}}_{n+1}$ if and only if
either $A  \in \overline{\mathcal{Z}}_{n+1}$ or $det(A)< 0,$
$ad(A)\geq 0$ and all principal $n \times n$ submatrices of $A$
are positive semi-definite.
\end{enumerate}
\end{theorem}

The paper is organized as follows. In section 2, we characterize
normal vectors of degenerated Euclidean simplexes. In section 3,
we characterize angle Gram matrices of degenerated hyperbolic
simplexes. Theorem 1 is proved in section 4 and Theorem 3 is
proved in section 5.

\subsection*{Acknowledgment}

We would like to thank the referee for a very careful reading of
the manuscript and for his/her nice suggestions on improving the
exposition.

\section{normal vectors of Euclidean simplexes}

As a convention, all vectors in $\mathbf{R}^m$ are column vectors
and the standard inner product in $\mathbf{R}^m$ is denoted by $u
\cdot v$. In the sequel, for a non-zero vector $w \in
\mathbf{R}^n$, we call the set $\{x\in \mathbf{R}^n | w\cdot x\geq
0 \}$ a \it closed half space\rm, and the set $\{x\in \mathbf{R}^n
| w\cdot x> 0 \}$ an \it open half space\rm. Define
$\mathcal{E}_{n+1}=\{(v_1,...,v_{n+1})\in (\mathbf{R}^n)^{n+1}|
v_1,...,v_{n+1}$ form unit outward normal vectors to the
codimension-1 faces of a Euclidean $n$-simplex $\}.$ Following
Milnor \cite{Mi}, a matrix is called \it unidiagonal \rm if its
diagonal entries are 1. An $(n+1)\times (n+1)$ symmetric
unidiagonal matrix $A$ is in $\mathcal{Z}_{n+1}$ if and only if
$A=[v_i\cdot v_j]$ for some point $(v_1,...,v_{n+1}) \in
\mathcal{E}_{n+1}$ (this is proved in \cite{Lu1},\cite{Mi}). We
claim that an $(n+1)\times (n+1)$ symmetric unidiagonal matrix $A$
is in $\overline{\mathcal{Z}}_{n+1}$ if and only if $A=[v_i\cdot
v_j]$ for some point $(v_1,...,v_{n+1})$ in the closure
$\overline{\mathcal{E}}_{n+1}$ in $(\mathbf{R}^n)^{n+1}.$ Indeed,
if $A=[v_i\cdot v_j]$ for some point $(v_1,...,v_{n+1}) \in
\overline{\mathcal{E}}_{n+1},$ then there is a sequence
$(v_1^{(m)},...,v_{n+1}^{(m)})\in \mathcal{E}_{n+1}$ converging to
$(v_1,...,v_{n+1}).$ We have a sequence of matrices
$A^{(m)}=[v_i^{(m)}\cdot v_j^{(m)}] \in \mathcal{Z}_{n+1}$
converging to $A.$ Conversely if $A \in
\overline{\mathcal{Z}}_{n+1}$, then there is a sequence of
matrices $A^{(m)}\in \mathcal{Z}_{n+1}$ converging to $A.$ Write
$A^{(m)}=[v_i^{(m)}\cdot v_j^{(m)}]$, where
$(v_1^{(m)},...,v_{n+1}^{(m)})\in \mathcal{E}_{n+1}$. Since
$v_i^{(m)}$ has norm 1 for all $i,m,$ by taking subsequence, we
may assume $\lim_{m\to
\infty}(v_1^{(m)},...,v_{n+1}^{(m)})=(v_1,...,v_{n+1})\in
\overline{\mathcal{E}}_{n+1}$ so that $A=[v_i\cdot v_j].$

A geometric characterization of elements in $\mathcal{E}_{n+1}$
was obtained in \cite{Lu1}. For completeness, we include a proof
here.

\begin{lemma} A collection of unit vectors
$(v_1,...,v_{n+1})\in (\mathbf{R}^n)^{n+1}$ is in
$\mathcal{E}_{n+1}$ if and only if one of the following conditions
is satisfied.

\begin{enumerate}
\item[(4.1)] The vectors $v_1,...,v_{n+1}$ are not in any closed
half-space.

\item[(4.2)] Any $n$ vectors of $v_1,...,v_{n+1}$ are linear
independent and the linear system $\sum_{i=1}^{n+1}a_iv_i=0$ has a
solution $(a_1,...,a_n)$ so that $a_i>0$ for all $i=1,...,n+1$.
\end{enumerate}
\end{lemma}

\begin{proof}
$(4.2)\Rightarrow(4.1).$ Suppose otherwise, $v_1,...,v_{n+1}$ are
in a closed half-space, i.e., there is a non-zero vector $w\in
\mathbf{R}^n$ so that $w\cdot v_i\geq 0, i=1,...,n+1.$  Let $a_1,
..., a_{n+1}$ be the positive numbers given by (4.2) so that
$\sum_{i=1}^{n+1} a_i v_i =0$. Then
$$0= w \cdot( \sum_{i=1}^{n+1} a_i v_i) = \sum_{i=1}^{n+1} a_i
(w \cdot v_i).$$ But by the assumption $a_i> 0,$ $w\cdot v_i\geq
0$ for all $i.$ Thus $w\cdot v_i=0$ for all $i.$ This means that
$v_1,...,v_{n+1}$ lie in the $(n-1)$-dimensional subspace
perpendicular to $w.$ It contradicts the assumption in (4.2) that
any $n$ vectors of $v_1,...,v_{n+1}$ are linear independent.

$(4.1)\Rightarrow(4.2).$ To see that any $n$ vectors of
$v_1,...,v_{n+1}$ are linear independent, suppose otherwise, some
$n$ vectors of $v_1,...,v_{n+1}$ are linear dependent. Therefore
there is an $(n-1)$-dimensional hyperplane containing these $n$
vectors. Then $v_1,...,v_{n+1}$ are contained in one of the two
closed half spaces bounded by the hyperplane. It contradicts to
the assumption of (4.1).

Since $v_1,...,v_{n+1}$ are linear dependent, and any $n$ of them
are linear independent, we can find real numbers $a_i\neq 0$ for
all $i$ such that $\sum_{i=1}^{n+1}a_iv_i=0.$ For any $i\neq j,$
let $H_{ij}$ be the $(n-1)$-dimensional hyperplane spanned by the
$n-1$ vectors $\{v_1,...,v_{n+1}\}\setminus\{v_i,v_j\}$ and $u\in
\mathbf{R}^n-\{0\}$ be a vector perpendicular to $H_{ij}.$ We have
$$0=u\cdot(\sum_{i=1}^{n+1}a_iv_i)=a_i(u\cdot v_i)+a_j(u\cdot
v_j).$$ By the assumption of (4.1), $v_i$ and $v_j$ must lie in
the different sides of $H_{ij}.$ Thus $u\cdot v_i$ and $u\cdot
v_j$ have different sign. This implies that $a_i$ and $a_j$ have
the same sign. Hence we can make $a_i>0$ for all $i$.

$\mathcal{E}_{n+1}\Leftrightarrow(4.1).$ We will show
$(v_1,...,v_{n+1})\in\mathcal{E}_{n+1}$ if and only if the
condition (4.1) holds. In fact, given an $n-$dimensional Euclidean
simplex $\sigma$, let $S^{n-1}$ be the sphere inscribed to
$\sigma.$ We may assume after a translation and a scaling that
$S^{n-1}$ is the unit sphere centered at the origin. Then the unit
vectors $v_1,...,v_{n+1}$ are the tangent points of $S^{n-1}$ to
the codimension-1 faces of $\sigma.$ The tangent planes to
$S^{n-1}$ at $v_i'$s bound a compact region (the Euclidean simplex
$\sigma$) containing the origin if and only if the tangent points
$v_1,...,v_{n+1}$ are not in any closed hemisphere of $S^{n-1}$.
\end{proof}

\begin{lemma} A collection of unit vectors
$(v_1,...,v_{n+1})\in (\mathbf{R}^n)^{n+1}$ is in
$\overline{\mathcal{E}}_{n+1}$ if and only if one of the following
conditions is satisfied:

\begin{enumerate}
\item [(5.1)] The vectors $v_1,...,v_{n+1}$ are not in any open
half-space.

\item [(5.2)] The linear system $\sum_{i=1}^{n+1}a_iv_i=0$ has a
nonzero solution $(a_1,...,a_{n+1})$ so that $a_i\geq 0$ for all
$i=1,...,n+1$.
\end{enumerate}
\end{lemma}

\begin{proof}
$\overline{\mathcal{E}}_{n+1}\Rightarrow$(5.1). To see that
elements in $\overline{\mathcal{E}}_{n+1}$ satisfy (5.1), if
$(v_1,...,v_{n+1}) \in \overline{\mathcal{E}}_{n+1},$ there is a
family of $(v_1^{(m)},...,v_{n+1}^{(m)})\in \mathcal{E}_{n+1}$
converging to $(v_1,...,v_{n+1}).$ Since vectors
$v_1^{(m)},...,v_{n+1}^{(m)}$ are not in any closed half-space for
any $m,$ by continuity, vectors $v_1,...,v_{n+1}$ are not in any
open half-space.

(5.1)$\Rightarrow$(5.2). Consider the linear map
$$f:\mathbf{R}^n\longrightarrow \mathbf{R}^{n+1}$$
$$w\longmapsto f(w)=[v_1,v_2,...,v_{n+1}]^t w=\left(
\begin{array}{cccc}
v_1\cdot w \\
v_2\cdot w \\
\vdots \\
v_{n+1}\cdot w
\end{array}
\right).$$ Statement (5.1) says that
\begin{eqnarray*}
\emptyset&=&\{w\in \mathbf{R}^n|v_i\cdot w>0, i=1,...,n+1\}\\
&=&\{w\in \mathbf{R}^n|f(w)>0\}\\
&=&f(\mathbf{R}^n)\cap \mathbf{R} ^{n+1}_{> 0}.
\end{eqnarray*}
Since $f(\mathbf{R}^n)$ and $\mathbf{R} ^{n+1}_{> 0}$ are convex
and disjoint, by the separation theorem for convex sets, there is
a vector $a=(a_1,...,a_{n+1})^t$ satisfying the conditions (i) and
(ii) below.

(i) For all $u\in \mathbf{R} ^{n+1}_{> 0},$
$$a\cdot u>0.$$
and

(ii) For all $w\in \mathbf{R} ^n,$ $$0\geq a\cdot f(w)= \left(
\begin{array}{cccc}
a_1\\
a_2\\
\vdots \\
a_{n+1}
\end{array}
\right) \cdot \left(
\begin{array}{cccc}
v_1\cdot w \\
v_2\cdot w \\
\vdots \\
v_{n+1}\cdot w
\end{array}
\right)=(\sum_{i=1}^{n+1}a_iv_i)\cdot w.
$$

Condition (i) implies that $a_i\geq 0,$ for $i=1,...,n+1$ and
$a\neq 0.$ Condition (ii) implies $\sum_{i=1}^{n+1}a_iv_i=0.$ Thus
(5.2) holds.

(5.2)$\Rightarrow\overline{\mathcal{E}}_{n+1}$. To see that a
point $(v_1,...,v_{n+1})$ satisfying (5.2) is in
$\overline{\mathcal{E}}_{n+1}$, we show that in any
$\varepsilon$-neighborhood of $(v_1,...,v_{n+1})$ in
$(\mathbf{R}^n)^{n+1}$, there is a point
$(v_1^\varepsilon,...,v_{n+1}^\varepsilon)\in \mathcal{E}_{n+1}.$

Let $\mathcal{N}_k$ be the set of $(v_1,...,v_k)$ such that
$v_i\in\mathbf{R}^{k-1},|v_i|=1$ for all $i$ and $\sum_{i=1}^k
a_iv_i=0$ has a nonzero solution $(a_1,...,a_k)$ with $a_i\geq 0$
for all $i$. The goal is to prove
$\mathcal{N}_{n+1}\subset\overline{\mathcal{E}}_{n+1}.$ We achieve
this by induction on $n.$ It is obvious that
$\mathcal{N}_2\subset\overline{\mathcal{E}}_2.$ Assume that
$\mathcal{N}_n\subset\overline{\mathcal{E}}_n$ holds.

For a point $(v_1,...,v_{n+1})\in\mathcal{N}_{n+1}$, if any $n$
vectors of $v_1,...,v_{n+1}$ are linear independent, then each
entry $a_i$ of the non-zero solution of the linear system
$\sum_{i=1}^{n+1}a_iv_i=0, a_i\geq 0, i=1,...,n+1$ must be
nonzero. Thus $a_i>0$ for all $i$ and $(v_1,...,v_{n+1})$
satisfies (4.2), therefore it is in $\mathcal{E}_{n+1}.$

In the remain case, without loss of generality, we assume that
$v_1,...,v_n$ are linear dependent. We may assume after a change
of coordinates that $v_i\in
\mathbf{R}^{n-1}=\mathbf{R}^{n-1}\times \{0\}\subset
\mathbf{R}^n,$ for $i=1,...,n,$ and
$v_{n+1}=(u_{n+1}\cos(\theta),\sin(\theta))^t,$ where
$0\leq\theta\leq\frac\pi2$ and $|u_{n+1}|=1.$

We claim that there exists some $1\leq i\leq n+1$ such that
$(v_1,...,\widehat{v}_i,$ $...,v_{n+1})\in\mathcal{N}_n,$ where
$\widehat{x}$ means deleting the element $x.$

Case 1. If $\theta>0$ i.e., $v_{n+1}$ is not in
$\mathbf{R}^{n-1},$ consider the nonzero solution of the linear
system $\sum_{i=1}^{n+1}a_iv_i=0, a_i\geq 0, i=1,...,n+1.$ The
last coordinate gives $a_10+...+a_n0+a_{n+1}\sin(\theta)=0,$ which
implies $a_{n+1}=0.$ This means $(a_1,...,a_n)\neq0,$ i.e.,
$(v_1,...,v_n)\in\mathcal{N}_n.$

Case 2. If $\theta=0$ i.e., $v_{n+1}\in \mathbf{R}^{n-1},$ then
the dimension of the solution space $W=\{(a_1,...,a_{n+1})^t\in
\mathbf{R}^{n+1}| \sum_{i=1}^{n+1}a_iv_i=0\}$ is at least 2. Since
$(v_1,...,v_{n+1})\in\mathcal{N}_{n+1},$ the intersection $W\cap
\mathbf{R}^{n+1}_{\geq0}-\{(0,...,0)\}$ is nonempty. The vector
space $W$ must intersect the boundary of the cone
$\mathbf{R}^{n+1}_{\geq0}-\{(0,...,0)\}$. Let $(a_1,...,a_{n+1})$
be a point in both $W$ and the boundary of the cone
$\mathbf{R}^{n+1}_{\geq0}-\{(0,...,0)\}$. Then there is some
$a_i=0.$ Then
$\{v_1(t),...,\widehat{v}_i(t),...,v_{n+1}(t)\}\in\mathcal{N}_n.$

By the above discussion, without lose of generality, we may assume
that $(v_1,...,v_n)\in\mathcal{N}_n$. By the induction assumption
$\mathcal{N}_n\subset\overline{\mathcal{E}}_n$, i.e., in any
$\frac\epsilon2$-neighborhood of $(v_1,...,v_n)$, we can find a
point $(u_1,...,u_n)\in \mathcal{E}_n,$ where
$u_i\in\mathbf{R}^{n-1}$ for all $i.$ Recall we write
$v_{n+1}=(u_{n+1}\cos(\theta),\sin(\theta))^t$. Let us define a
continuous family of $n+1$ unit vectors $v_1(t),...,v_{n+1}(t)$ by
setting
$$v_i(t)=(u_i\cos(t^2),-\sin(t^2))^t,1\leq i\leq n,$$
$$v_{n+1}(t)=(u_{n+1}\cos(\theta+t),\sin(\theta+t))^t.$$
We claim that there is a point $(v_1(t),...,v_{n+1}(t))\in
\mathcal{E}_{n+1}$ for small $t>0$ within
$\frac\epsilon2$-neighborhood of
$((u_1,0)^t,...,(u_n,0)^t,v_{n+1}).$ By triangle inequality, this
point is within $\varepsilon$-neighborhood of $(v_1,...,v_{n+1})$.
We only need to check that $(v_1(t),...,v_{n+1}(t))\in
\mathcal{E}_{n+1}$ for sufficiently small $t>0$ by verifying the
condition (4.2).

To show any $n$ vectors of $v_1(t),...,v_{n+1}(t)$ are linear
independent, it is equivalent to show that
$$\det[v_1(t),...,\widehat{v}_i(t),...,v_{n+1}(t)]\neq0$$ for each
$i=1,...,n+1.$

First,
\begin{eqnarray*}
\det[v_1(t),...,v_n(t)]&=&\det \left[
\begin{array}{ccccc}
u_1\cos(t^2) & u_2\cos(t^2) & ... & u_n\cos(t^2) \\
-\sin(t^2) & -\sin(t^2) & ... & -\sin(t^2)\\
\end{array}
\right]_{n\times n} \\
&=&-\sin(t^2)\cos(t^2)^{n-1}\det \left[
\begin{array}{ccccc}
u_1 & u_2 & ... & u_n \\
1 & 1 & ... & 1 \\
\end{array}
\right].
\end{eqnarray*}

To see the determinant is nonzero, suppose there are real numbers
$a_1,...,a_n$ such that $\sum_{i=1}^na_i(u_1,1)^t=0.$ Then we have
$\sum_{i=1}^na_iu_i=0$ and $\sum_{i=1}^na_i=0$. By assumption
$(u_1,...,u_n) \in \mathcal{E}_n$, we know either $a_i=0$ for all
$i$ or $a_i\neq 0$ and have the same sign for all $i$. Hence
$\sum_{i=1}^na_i=0$ implies $a_i=0$ for all $i$. Thus the vectors
$(u_1,1)^t,...,(u_n,1)^t$ are linear independent. Hence
$\det[v_1(t),...,v_n(t)]\neq 0$ for $t\in (0, \sqrt{\frac\pi2} ].$

Second, we calculate the determinant of the matrix whose columns
are $v_{n+1}(t)$ and some $n-1$ vectors of $v_1(t),...,v_n(t).$
Without loss of generality, consider
$$f(t)=\det[v_2(t),...,v_n(t),v_{n+1}(t)]$$
$$=\det \left[
\begin{array}{ccccc}
u_2\cos(t^2) & ... & u_n\cos(t^2) & u_{n+1}\cos(\theta+t)\\
-\sin(t^2) & ... & -\sin(t^2) & \sin(\theta+t)\\
\end{array}
\right] $$ If $\theta>0,$ by the assumption $(u_1,...,u_n) \in
\mathcal{E}_n$, then
$$
f(0)=\det \left[
\begin{array}{ccccc}
u_2& ... & u_n & u_{n+1}\cos(\theta) \\
0 & ... & 0 & \sin(\theta) \\
\end{array}
\right] =\sin(\theta)\det[u_2,...,u_n]\neq0.
$$
It implies $f(t)\neq 0$ holds for sufficiently small $t>0.$

If $\theta=0,$ then $f(0)=0.$ By expanding the determinant,
$$f(t)=-\sin(t^2)g(t)+\sin(t)\det[u_2\cos(t^2),...,u_n\cos(t^2)],$$ for
some function $g(t)$, therefore $f'(0)=\det(u_2,...,u_n)\neq0.$
Hence $f(t)\neq 0$ holds for sufficiently small $t>0.$

Next, let $$a_i(t)=(-1)^{i-1}\det[v_1(t),...,\widehat{v}_i(t),
...,v_{n+1}(t)],1\leq i\leq n+1,$$ then
$\sum_{i=1}^{n+1}a_i(t)v_i(t)=0.$ Since
$\det[v_1(0),...,v_n(0)]=0,$ we have $a_{n+1}(0)=0.$ This shows
that $\sum_{i=1}^na_i(0)v_i(0)=0,$ therefore
$\sum_{i=1}^na_i(0)u_i=0$. By the assumption $(u_1,...,u_n) \in
\mathcal{E}_n$, we obtain $a_i(0)\cdot a_j(0)>0$ for $0\leq
i,j\leq n.$ By the continuity we obtain $a_i(t)\cdot a_j(t)>0$ for
$0\leq i,j\leq n,$ for sufficient small $t>0.$ Consider the last
coordinate of $\sum_{i=1}^{n+1}a_i(t)v_i(t)=0$ we obtain
$$-\sin(t^2)\sum_{i=1}^na_i(t)+\sin(\theta+t)a_{n+1}(t)=0.$$ Thus
$a_{n+1}(t)$ has the same sign as that of $a_i(t).$ Thus
$(a_1(t),...,a_{n+1}(t))$ or $(-a_1(t),...,-a_{n+1}(t))$ is a
solution required in condition (4.2).
\end{proof}

\section{degenerate hyperbolic simplexes}

Let $\mathbf{R}^{n,1}$ be the Minkowski space which is
$\mathbf{R}^{n+1}$ with an inner product $\langle,\rangle$ where
$$\langle(x_1,...,x_n,x_{n+1})^t,(y_1,...,y_n,y_{n+1})^t\rangle
=x_1y_1+...+x_ny_n-x_{n+1}y_{n+1}.$$ Let
$H^n=\{x=(x_1,...,x_{n+1})^t\in \mathbf{R}^{n,1}|\langle
x,x\rangle=-1,x_{n+1}>0\}$ be the hyperboloid model of the
hyperbolic space. The de Sitter space is
$\{x\in\mathbf{R}^{n,1}|\langle x,x\rangle=1\}.$ For a hyperbolic
simplex $\sigma$ in $H^n,$ the \it center \rm and the \it radius
\rm of the simplex $\sigma$ are defined to be the center and
radius of its inscribed ball.

\begin{lemma} For an n-dimensional hyperbolic simplex $\sigma\in H^n$
with center $e_{n+1}=(0,...,0,1)^t$, its unit outward normal
vectors in the de Sitter space are in a compact set independent of
$\sigma.$
\end{lemma}

\begin{proof} Let $v_1,...,v_{n+1}$ be the unit outward normal vectors of
$\sigma,$ i.e., $$\sigma=\{x\in H^n|\langle x,v_i\rangle\leq0\
\mbox{and}\ \langle v_i,v_i\rangle=1\ \mbox{for all}\ i\}.$$ Let
$v_i^\bot$ be the totally geodesic hyperplane in $H^n$ containing
the $(n-1)-$dimensional face of $\sigma$ perpendicular to $v_i$
for each $i=1,...,n+1.$ The radius of $\sigma$ is the distance
from the center $e_{n+1}$ to $v_i^\bot$ for any $i=1,...,n+1$
which is equal to $\sinh^{-1}(|\langle e_{n+1},v_i\rangle|)$ (see
for instance \cite{Vi2} p26). It is well known that the volume of
an $n-$dimensional hyperbolic simplex is bounded by the volume of
the $n-$dimensional regular ideal hyperbolic simplex which is
finite (see for instance \cite{R} p539). It implies that the
radius of a hyperbolic simplex $\sigma$ is bounded from above by a
constant independent of $\sigma.$ Hence $\langle
e_{n+1},v_i\rangle^2$ is bounded from above by a constant $c_n$
independent of $\sigma$ for any $i=1,...,n+1$. It follows that
$v_1,...,v_{n+1}$ are in the compact set
\begin{eqnarray*}
X_n &=& \{x=(x_1,...,x_{n+1})^t| \langle x,x\rangle=1,\langle
e_{n+1},x\rangle^2\leq c_n\}\\
&=&\{x=(x_1,...,x_{n+1})^t| x_1^2+...+x_n^2=x_{n+1}^2+1,
x_{n+1}^2\leq c_n\}
\end{eqnarray*} independent of $\sigma.$
\end{proof}

\begin{lemma} If $A \in
\overline{\mathcal{Y}}_{n+1}$ and $\det(A)=0$, then $A \in
\overline{\mathcal{Z}}_{n+1}$.
\end{lemma}

\begin{proof} Let $A^{(m)}$ be a sequence of angle Gram matrixes in
$\mathcal{Y}_{n+1}$ converging to $A.$ By Proposition 2 (c), for
any $m$, all principal $n\times n$ submatrices of $A^{(m)}$ are
positive definite. Thus all principal $n\times n$ submatrices of
$A$ are positive semi-definite. Since $\det(A)=0,$ we see that $A$
is positive semi-definite.

Let $\sigma^{(m)}$ be the $n-$dimensional hyperbolic simplex in
the hyperboloid model $H^n$ whose angle Gram matrix is $A^{(m)}$
and whose center is $e_{n+1}=(0,...,0,1)^t$. By Lemma 6, its unit
outward normal vector $v_i^{(m)}$ is in a compact set. Thus by
taking a subsequence, we may assume
$(v_1^{(m)},...,v_{n+1}^{(m)})$ converges to $(v_1,...,v_{n+1})$
with $\langle v_i,v_i\rangle=1$. Since $A^{(m)}=[\langle
v_i^{(m)},v_j^{(m)}\rangle]$ and $A^{(m)}$ converges to $A,$ we
obtain $$A=[\langle
v_i,v_j\rangle]=[v_1,...,v_{n+1}]^tS[v_1,...,v_{n+1}],$$ where $S$
is the diagonal matrix $diag(1,...,1,-1).$

Since $\det(A)=0,$ the vectors $v_1,...,v_{n+1}$ are linear
dependent. Assume that the vectors $v_1,...,v_{n+1}$ span a
$k-$dimensional subspace $W$ of $\mathbf{R}^{n,1}$, where $k\leq
n.$

For any vector $x\in W,$ write $x=\sum_{i=1}^{n+1}x_iv_i.$ Then
\begin{eqnarray*}
\langle x,x
\rangle&=&(x_1,...,x_{n+1})[v_1,...,v_{n+1}]^tS[v_1,...,v_{n+1}](x_1,...,x_{n+1})^t\\
&=&(x_1,...,x_{n+1})A(x_1,...,x_{n+1})^t\\
&\geq&0
\end{eqnarray*}
due to the fact that $A$ is positive semi-definite.

Now for any $x,y\in W,$ the inequality $\langle
x+ty,x+ty\rangle\geq0$ for any $t\in \mathbf{R}$ implies the
Schwartz inequality
$$\langle x,y\rangle^2\leq\langle x,x\rangle\langle y,y\rangle.$$

To prove that $A \in \overline{\mathcal{Z}}_{n+1}$, we consider
the following two possibilities.

Case 1. If $\langle x,x \rangle>0$ holds for any non-zero $x\in
W,$ then the Minkowski inner product restricted on $W$ is positive
definite. Since the Minkowski inner product restricted on
$\mathbf{R}^k=\mathbf{R}^k\times 0\subset \mathbf{R}^{n,1}$ is
positive definite, by Witt's theorem, there is an isometry
$\gamma$ of $\mathbf{R}^{n,1}$ sending $W$ to $\mathbf{R}^k$ (see
\cite{Vi2} p14-p15). By replacing $v_i^{(m)}$ by
$\gamma(v_i^{(m)})$ for each $i$ and $m$, we may assume that
$v_1,...,v_{n+1}$ are contained in $\mathbf{R}^k$. Thus $\langle
v_i,v_i\rangle=v_i\cdot v_j$ for all $i,j.$ Therefore
$$A=[v_1,...,v_{n+1}]^tS[v_1,...,v_{n+1}]=[v_1,...,v_{n+1}]^t[v_1,...,v_{n+1}].$$

To show $A\in\overline{\mathcal{Z}}_{n+1},$ by Lemma 5, we only
need to show that $v_1,...,v_{n+1}$ are not contained in any open
half space of $\mathbf{R}^k.$ This is the same as that
$v_1,...,v_{n+1}$ are not contained in any open half space of
$\mathbf{R}^n$.

Suppose otherwise, there exists a vector $w\in\mathbf{R}^k$ such
that $v_i\cdot w>0$ for all $i.$ Thus $\langle v_i,
w\rangle=v_i\cdot w>0.$ By taking $m$ large enough, we obtain
$\langle v_i^{(m)}, w\rangle>0$ for all $i.$

It is well known that for the unit normal vectors $v_i^{(m)}$ of a
compact hyperbolic simplex in $H^n$, the conditions  $\langle
v_i^{(m)}, w\rangle>0$ for all $i$ implies $\langle w, w\rangle
<0.$ But this contradicts the assumption that $w\in\mathbf{R}^k$
which implies $\langle w, w\rangle\geq 0.$

Case 2. If there exists some non-zero vector $x_0\in W$ such that
$\langle x_0,x_0 \rangle=0,$ then by the Schwartz inequality we
have $$\langle x_0,y\rangle^2\leq\langle x_0,x_0\rangle\langle
y,y\rangle=0$$ for any $y\in W.$ Thus $\langle x_0,y\rangle=0$ for
any $y\in W.$ This implies that the subspace $W$ is contained in
$x_0^\perp$, the orthogonal complement of $x_0.$

Since the vector $u=(0,...,0,1,1)^t\in\mathbf{R}^{n,1}$ satisfies
$\langle u,u\rangle=0,$ there is an isometry $\gamma$ of
$\mathbf{R}^{n,1}$ sending $x_0$ to $u.$ Thus $\gamma$ sends
$x_0^\perp$ to $u^\perp.$ By replacing $v_i^{(m)}$ by
$\gamma(v_i^{(m)})$ for each $i$ and $m$, we may assume that
$v_1,...,v_{n+1}$ are contained in $u^\perp.$

For any $i,$ since $\langle v_i,u\rangle=\langle
v_i,(0,...,0,1,1)^t\rangle=0$, we can write $v_i$ as
$$v_i=w_i+a_iu$$ for some $w_i\in\mathbf{R}^{n-1}$ and
$a_i\in\mathbf{R}.$ Since $\langle w_i,u\rangle=0,$ thus $\langle
v_i,v_j\rangle=w_i\cdot w_j$ for all $i,j.$ Therefore
\begin{eqnarray*}
A&=&[v_1,...,v_{n+1}]^tS[v_1,...,v_{n+1}]\\
&=&[w_1,...,w_{n+1}]^t[w_1,...,w_{n+1}].
\end{eqnarray*}

To show $A\in\overline{\mathcal{Z}}_{n+1},$ by Lemma 5, we only
need to show that $w_1,...,w_{n+1}$ are not contained in any open
half space of $\mathbf{R}^{n-1}$ which is equivalent to that
$w_1,...,w_{n+1}$ are not contained in any open half space of
$\mathbf{R}^n$.

Suppose otherwise, there exists a vector $w\in\mathbf{R}^{n-1}$
such that $w_i\cdot w>0$ for all $i.$ Then $$\langle v_i,
w\rangle=\langle w_i, w\rangle+\langle (0,...,0,a_i,a_i)^t,
w\rangle=w_i\cdot w+0>0$$ for all $i.$ By taking $m$ large enough,
we obtain $\langle v_i^{(m)}, w\rangle>0$ for all $i.$ By the same
argument above, it is a contradiction.

\end{proof}

\section{proof of theorem 1}

\subsection*{Spherical case}
We begin with a brief review of the relevant result in \cite{Lu2}.
For any positive semi-definite symmetric matrix $A$, there exists
a unique positive semi-definite symmetric matrix $\sqrt{A}$ so
that $(\sqrt{A})^2 = A.$ It is well know that the map $A
\longmapsto \sqrt{A}$ is continuous on the space of all positive
semi-definite symmetric matrices.

Suppose $A \in \mathcal{X}_{n+1} =\{ A=[a_{ij}] \in\mathbf{R}^{
(n+1) \times (n+1)} |$
 $A^t = A$, all $a_{ii} =1$, $A$ is positive definite\}, the space of the angle Gram
matrices of spherical simplexes (by the Proposition 2). By making
a change of variables, the Aomoto-Kneser-Vinberg formula (\ref{1}) is
equivalent to
\begin{equation} \label{2}
V(A) = \mu_n^{-1} \int_ {\mathbf{R} ^{n+1}}\chi(\sqrt{A}x)e^{ -x^t x}dx,
\end{equation}
where $\chi$ is the characteristic function of the set $
\mathbf{R}^{n+1}_{\geq 0}$ in $ \mathbf{ R}^{n+1}.$ It is proved
in \cite{Lu2} that volume formula (\ref{2}) still holds for any
matrix in $\overline{\mathcal{X}}_{n+1}$=\{$A=[a_{ij}]
\in\mathbf{R}^{ (n+1) \times (n+1)} |$ $A^t = A$, all $a_{ii} =1$,
$A$ is positive semi-definite\}.

Suppose $V(A)=0,$ by formula (\ref{2}), we see the function
$\chi\circ h:\mathbf{R}^{n+1}\to\mathbf{R}$ is zero almost
everywhere, where $h:\mathbf{R}^{n+1}\to\mathbf{R}^{n+1}$ is the
linear map sending $x$ to $\sqrt{A}x.$ Equivalently, the
$(n+1)$-dimensional Lebesque measure of
$h^{-1}(\mathbf{R}^{n+1}_{\geq 0})$ is zero. We claim
$h(\mathbf{R}^{n+1})\cap \mathbf{R}^{n+1}_{> 0}=\emptyset.$ For
otherwise, $h^{-1}(\mathbf{R} ^{n+1}_{> 0})$ is a nonempty open
subset in $\mathbf{R}^{n+1}$ with positive $(n+1)$-dimensional
Lebesque measure. This is a contradiction.

Now let $\sqrt{A}=[v_1,...,v_{n+1}]^t_{(n+1)\times(n+1)},$ where
$v_i\in \mathbf{R} ^{n+1}$ is a column vector for each $i$. First
$h(\mathbf{R}^{n+1})\cap \mathbf{R}^{n+1}_{> 0}=\emptyset$ implies
that $\det\sqrt{A}=0.$ Therefore $\{v_1,...,v_{n+1}\}$ are linear
dependent. We may assume, after a rotation $r\in O(n+1)$, the
vectors $v_1,...,v_{n+1}$ lie in $\mathbf{R} ^n\times \{0\}.$ Now
\begin{eqnarray*} \emptyset&=&h(\mathbf{R}^{n+1})\cap \mathbf{R}^{n+1}_{>0}\\
&=&\{\sqrt{A}w|w\in\mathbf{R}^{n+1}\}\cap \mathbf{R}^{n+1}_{>0}\\
&=&\{(v_1\cdot w,...,v_{n+1}\cdot w)^t|w\in\mathbf{R}^{n+1}\}\cap \mathbf{R}^{n+1}_{>0}
\end{eqnarray*}
This shows that there is no $w\in\mathbf{R}^{n+1}$ such that
$v_i\cdot w>0$ for $i=1,...,n+1,$ i.e., the vectors
$v_1,...,v_{n+1}$ are not in any open half space. By lemma 5, we
have $(v_1,...,v_{n+1}) \in \overline{\mathcal{E}}_{n+1},$
therefore $A=[v_i\cdot v_j]\in \overline{\mathcal{Z}}_{n+1}.$

\subsection*{Hyperbolic case}
Let $A\in \overline{\mathcal{Y}}_{n+1}$. If $\det(A)\neq0,$ it is
proved in \cite{Lu2} that the volume formula (\ref{1}) still holds
for $A.$ In formula (\ref{1}), since $-x^t ad(A) x$ is finite, the
integrant $e^{ -x^t ad(A) x}> 0$. Therefore the integral
$\int_{\mathbf{R} ^{n+1}_{\geq 0}} e^{ -x^t ad(A) x} dx > 0$.
Hence $V(A)>0.$

It follows that if the extended volume function vanishes at $A$,
then $\det{A}=0.$ By Lemma 7 , we have $A\in
\overline{\mathcal{Z}}_{n+1}.$

\section{proof of theorem 3}

\subsection*{Proof of (a)}
If $A \in \overline{\mathcal{Z}}_{n+1},$ then $A=[v_i\cdot v_j]$
for some point $(v_1,...,v_{n+1}) \in
\overline{\mathcal{E}}_{n+1}.$ By Lemma 5, the linear system
$\sum_{i=1}^{n+1}a_iv_i=0, a_i\geq 0, i=1,...,n+1$ has a nonzero
solution. Let $(a_1,...,a_{n+1})$ be a solution with the least
number of nonzero entries among all solutions. By rearrange the
index, we may assume $a_1>0,...,a_{k+1}>0, a_{k+2}=...=a_{n+1}=0.$
We claim $rank[v_1,...,v_{k+1}]=k.$ Otherwise
$rank[v_1,...,v_{k+1}]\leq k-1,$ then the dimension of the
solution space $W=\{(x_1,...,x_{k+1})^t\in \mathbf{R}^{k+1}|
\sum_{i=1}^{k+1}x_iv_i=0\}$ is at least 2. Thus $\Omega=W\cap
  \mathbf{R}^{k+1}_{>0}$ is
a nonempty open convex set in $W$ whose dimension is at least 2.
Hence $\Omega$ contains a boundary point $(b_1,...,b_{k+1})\in
\Omega-\{(0,...,0)\}$ with some $b_j=0,$ due to $dim W\geq 2.$ Now
we obtain a solution
$(b_1,...b_{j-1},0,b_{j+1},...,b_{k+1},0,...,0)$ which has lesser
number of nonzero entries than $(a_1,...,a_{n+1}).$ This is a
contradiction.

Let $B=[v_i\cdot v_j]_{(k+1)\times(k+1)}.$ Since
$rank[v_1,...,v_{k+1}]=k$, we have $\det(B)=0.$ We claim that
$ad(B)\geq 0$ and $ad(B)\neq 0.$ This will verify the condition
(a) in Theorem 3 for $A.$ Let $ad(B)=[b_{ij}]_{(k+1)\times(k+1)}.$
Evidently, due to $rank(B)=k, ad(B)\neq0.$ It remains to prove
that $ad(B)\geq0.$ By the construction of $B$, we see $b_{jj}\geq
0,$ for all $j.$ Since $rank[v_1,...,v_{k+1}]=k$, it follows the
dimension of the solution space of $\sum_{i=1}^{k+1}a_iv_i=0$ is
1. Since $\sum_{i=1}^{k+1}b_{ij}v_j=0,$ $(b_{i1},...,b_{ik+1})$ is
proportional to $(a_1,...,a_{k+1}),$ where $a_i>0$ for $1\leq
i\leq k+1.$ This shows that if $b_{jj}> 0,$ then $b_{ij}>0$ for
all $i$, if $b_{jj}= 0,$ then $b_{ij}= 0$ for all $i$. This shows
$ad(B)\geq0.$

Conversely, if $A$ is positive semi-definite so that $\det(A)=0$
and there exists a principal $(k+1)\times (k+1)$ submatrix $B$ so
that $\det(B)=0, ad(B)\geq 0$ and $ad(B)\neq 0$, we will show that
$A \in  \overline{\mathcal{Z}}_{n+1}$. Since $A$ is positive
semi-definite and unidiagonal, there exist unit vectors
$v_1,...,v_{n+1}$ in $\mathbf{R}^n$ such that $A= [v_i\cdot v_j].$
We may assume $B=[v_i\cdot v_j]_{(k+1)\times(k+1)}, 1\leq i,j \leq
k+1$ and $ad(B)=[b_{ij}].$ Due to $\det(B)=0, ad(B)\neq 0,$ we
have $rank(v_1,...,v_{k+1})=k.$ We may assume $v_2,...,v_{k+1}$
are independent. Thus the cofactor $b_{11}>0.$ By the assumption
$ad(B)\geq 0,$ we have $b_{1s}\geq 0$ for $s=1,...k+1.$ Since
$\sum_{s=1}^{k+1}b_{1s}(v_s\cdot v_j)=0$ for all $j=2,...,k+1$ and
$v_2,...,v_{k+1}$ are independent, we get
$\sum_{s=1}^{k+1}b_{1s}v_s=0.$ Thus we get a nonzero solution for
the linear system $\sum_{i=1}^{n+1}a_iv_i=0, a_i\geq 0,
i=1,...,n+1.$

\subsection*{Proof of (b)}
If $A\in \overline{\mathcal{X}}_{n+1}-\mathcal{X}_{n+1},$ then
$A=[v_i\cdot v_j]$ where $v_1,...,v_{n+1}$ are linear dependent.
We can assume $v_1,...,v_{n+1}$ lie in $\mathbf{R}^n\times\{0\}.$
By change subindex, we may assume $\sum_{i=1}^{n+1}a_iv_i=0$ has a
non-zero solution with $a_i\geq0$ if $i=1,...,k$ while $a_i<0$ if
$i=k+1,...,n+1.$ Thus vectors $v_1,...,v_k,-v_{k+1},...,-v_{n+1}$
satisfy the condition (5.2) in Lemma 5. Let $D$ be the diagonal
matrix $diag(1,...,1,-1,...,-1)$ with $k$ diagonal entries being 1
and $n-k+1$ diagonal entries being $-1.$ Thus by Lemma 5,
$DAD\in\overline{\mathcal{Z}}_{n+1}.$

On the other hand, if for some diagonal matrix $D$ in Theorem 3
(b), we have $DAD\in\overline{\mathcal{Z}}_{n+1},$ then by Theorem
3 (a), $DAD$ is positive semi-definite. Therefore $A$ is positive
semi-definite. Take $B\in \mathcal{X}_{n+1}$ and consider the
family $A(t)=(1-t)A+tB$ for $t\in [0,1].$ Then $\lim_{t \to
0}A(t)=A$ and $A(t)\in \mathcal{X}_{n+1}$ for $t>0.$ Thus $A\in
\overline{\mathcal{X}}_{n+1}.$

\subsection*{Proof of (c)}
First we show that the conditions are sufficient. Suppose
$A=[a_{ij}]_{(n+1)\times(n+1)}$ is a symmetric unidiagonal matrix
with all principal $n \times n$ submatrices positive semi-definite
so that either $A\in \overline{\mathcal{Z}}_{n+1}$ or $\det(A)< 0$
and $ad(A)\geq 0.$ We will show $A\in
\overline{\mathcal{Y}}_{n+1}.$ If $A\in
\overline{\mathcal{Z}}_{n+1},$ it is sufficient to show that $
\mathcal{Z}_{n+1}\subset \overline{\mathcal{Y}}_{n+1},$ i.e., we
may assume $A\in \mathcal{Z}_{n+1}.$ In this case, let
$J=[c_{ij}]_{(n+1)\times(n+1)}$ so that $c_{ii}=1$ and $c_{ij}=-1$
for $i\neq j.$ Consider the family $A(t)=(1-t)A+tJ,$ for $0\leq
t\leq 1.$ Evidently $\lim_{t\to 0}A(t)=A.$ We claim that $A(t)\in
\mathcal{Y}_{n+1}$ for small $t>0.$ Since all principal $n \times
n$ submatrices of $A$ are positive definite, by continuity, all
principal $n \times n$ submatrices of $A(t)$ are positive definite
for small $t>0.$ It remains to check $\det(A(t))<0$ for small
$t>0.$ To this end, let us consider
$\frac{d}{dt}|_{t=0}\det(A(t)).$ We have
$$\frac{d}{dt}|_{t=0}\det(A(t))=\sum_{i\neq j}(-a_{ij}-1)cof(A)_{ij}<0,$$
due to $ad(A)=[cof(A)_{ij}]>0$ and
$-a_{ij}-1<0$ for all $i\neq j.$ Since $\det(A)=0$  it follows that
$det(A(t))<0$ for small $t>0.$

In the second case that $\det(A)< 0$ and $ad(A)\geq 0$ and all principal $n
\times n$ submatrices of $A$ are positive semi-definite. Then $A$
has a unique negative eigenvalue $-\lambda$, where $\lambda>0.$
Consider the family $A(t)=A+t\lambda I,$ for $0\leq t\leq 1,$
where $I$ is the identity matrix, so that $$\lim_{t\to
0}\frac1{1+\lambda t}A(t)=A.$$ We claim there is a diagonal matrix
$D$ whose diagonal entries are $\pm 1$ so that

\begin{enumerate}
\item [(1)] $DAD=A$,

\item [(2)] $\frac1{1+\lambda t}DA(t)D\in \mathcal{Y}_{n+1}$ for
$0<t<1.$
\end{enumerate}

As a consequence, it follows
\begin{eqnarray*}
A&=&DAD\\
&=&\lim_{t\to 0}\frac1{1+\lambda t}DA(t)D
\in\overline{ \mathcal{Y}}_{n+1}.
\end{eqnarray*}
To find this diagonal
matrix $D,$ by the continuity, $\det(A(t))<0$ for $0<t<1$ and
$\det(A(1))=0.$ Furthermore, all principal $n \times n$ submatrices
of $A$ are positive definite for $t>0$ due to positive definiteness of $t\lambda I.$  Let us recall the Lemma
3.4 in \cite{Lu2} which says that if $B$ is a symmetric $(n+1) \times
(n+1)$ matrix so that all $n \times n$ principal submatrices in
$B$ are positive definite and $\det(B) \leq 0$, then no entry in
the adjacent matrix $ad(B)$ is zero. It follows that every entry
of $ad(A(t))$ is nonzero for $0<t\leq1.$

Let $ad(A(1))=[b_{ij}]_{(n+1)\times(n+1)}$ and $D$ to be the
diagonal matrix with diagonal entries being
$$\frac{b_{1i}}{|b_{1i}|}= \pm 1$$ for $i=1,...,n+1.$ Then the
entries of the first row and the first column of $Dad(A(1))D$ are
positive. Since $\det(A(1))=0$ and $ad(A(1))\neq 0$, we see the
rank of $ad(A(1))$ is 1. Thus any other column is propositional to
the first column. But $b_{ii}>0$ for all $i$, hence $ad(A(1))>0.$
Now since every entry of $Dad(A(t))D$ is nonzero for $t>0,$ by
continuity $Dad(A(t))D>0$ for $t>0$ and $Dad(A)D=Dad(A(0))D\geq
0.$ By the assumption $ad(A)\geq 0,$ it follows $Dad(A)D=ad(A).$
On the other hand $Dad(A)D=ad(D^{-1}AD^{-1})$, and $det(A)\neq 0.$
Thus $D^{-1}AD^{-1}=A$ or the same $A=DAD.$ This shows
\begin{eqnarray*}
A&=&DAD\\
&=&\lim_{t\to 0}DA(t)D\\
&=&\lim_{t\to 0}\frac1{1+\lambda t}DA(t)D.
\end{eqnarray*}
By the construction above $\frac1{1+\lambda t}DA(t)D\in
\mathcal{Y}_{n+1}$ for $0<t<1,$ this shows $A\in
\overline{\mathcal{Y}}_{n+1}.$

Finally, we show the condition in (c) is necessary. Suppose
$A=\lim_{m\to \infty}A^{(m)}$ where $A^{(m)}\in
\mathcal{Y}_{n+1}.$ By Proposition 2, $\det(A^{(m)})<0,$ $
ad(A^{(m)})>0$ and all principal $n\times n$ submatrices of
$A^{(m)}$ are positive definite. We want to show that $A$
satisfies the conditions stated in (c). Evidently, all principal
$n \times n$ submatrices of $A$ are positive semi-definite,
$ad(A)\geq 0$ and $det(A)\leq 0.$ If $\det(A)< 0,$ then we are
done. If $\det(A)= 0,$ by Lemma 7, we see that $A\in
\overline{\mathcal{Z}}_{n+1}$.

\bibliographystyle{amsplain}

\end{document}